\documentclass[11pt, a4paper, reqno]{amsart}
\usepackage{amssymb, array, amsmath, amscd, pdfpages, enumerate, amsthm, fixltx2e, setspace, hyperref,mathtools}
\pdfoutput=1
\usepackage[latin1]{inputenc} 
\usepackage{amssymb}
\usepackage{bm}
\usepackage{graphicx}

\DeclareMathOperator{\Span}{span}

\newcommand*{\R}{{\mathbb{R}}}

\newcommand*{\N}{{\mathbb{N}}}

\newcommand*{\Abs}[2][default]{\ifthenelse{\equal{#1}{default}}{\left\lvert#2\right\rvert}{\ldelim{#1}{\lvert}#2\rdelim{#1}{\rvert}}}
\newcommand*{\Norm}[2][default]{\ifthenelse{\equal{#1}{default}}{\left\lVert#2\right\rVert}{\ldelim{#1}{\lVert}#2\rdelim{#1}{\rVert}}}

\newcommand*{\Iprod}[3][default]{\ifthenelse{\equal{#1}{default}}{\left\langle#2,#3\right\rangle}{\ldelim{#1}{\langle}#2,#3\rdelim{#1}{\rangle}}}
\newcommand*{\Dualpair}[3][default]{\ifthenelse{\equal{#1}{default}}{\left\langle#2,#3\right\rangle}{\ldelim{#1}{\langle}#2,#3\rdelim{#1}{\rangle}}}

\newcommand*{\ddb}[2][1]{\ifthenelse{\equal{#1}{1}}{\frac{d}{d#2}}{\frac{d^{#1}}{d#2^{#1}}}}
\newcommand*{\pd}[3][1]{\ifthenelse{\equal{#1}{1}}{\frac{\partial{#2}}{\partial{#3}}}{\frac{\partial^{#1}{#2}}{\partial#3^{#1}}}}

\newcommand*\lenv{{\hbox{\raisebox{-.15ex}{\rotatebox[origin=c]{50}{$\smallsmile$}}\kern-8.65pt\rotatebox[origin=c]{-25}{$\smallsetminus$}}}}
\newcommand*\uenv{{\hbox{\raisebox{-.0ex}{\rotatebox[origin=c]{-45}{$\smallfrown$}}\kern-5.6pt\raisebox{.2ex}{\rotatebox[origin=c]{-5}{\scriptsize\slash}}}}\,\kern+1.5pt}

\newcommand*{\lupp}[1]{\kern0pt^{\kern0pt u \kern0pt}\kern1pt#1}
\newcommand*{\llow}[1]{\kern0pt^{\kern0pt l \kern0pt}\kern1pt#1}
\newcommand*{\rupp}[1]{#1\kern1pt^{\kern0pt u \kern0pt}\kern0pt}
\newcommand*{\rlow}[1]{#1\kern1pt^{\kern0pt l \kern0pt}\kern0pt}
\newcommand*{\ul}[1]{\kern0pt^{\kern0pt u \kern0pt}\kern0pt#1\kern1pt^{\kern0pt l \kern0pt}\kern0pt}
\newcommand*{\lu}[1]{\kern0pt^{\kern0pt l \kern0pt}\kern0pt#1\kern0pt^{\kern0pt u \kern0pt}\kern0pt}

\DeclareMathOperator{\mmin}{Min}
\DeclareMathOperator{\mmax}{Max}
\newcommand{\mini}[1]{\mmin [\kern1pt #1 \kern1pt ]}
\newcommand{\maxi}[1]{\mmax [\kern1pt #1 \kern1pt ]}

\newcommand{\pperp}{\perp^*}

\newcommand*{\bproofname}{Proof}
\newenvironment{bproof}[1][\bproofname]{\begin{proof}[#1]}{\end{proof}}

\newtheorem{thm}{Theorem}[section]
\newtheorem{prop}[thm]{Proposition}
\newtheorem{lemma}[thm]{Lemma}
\newtheorem{cor}[thm]{Corollary}
\theoremstyle{definition}
\newtheorem{defn}[thm]{Definition}

\newtheorem{remark}[thm]{Remark}
\newtheorem{example}[thm]{Example}

\numberwithin{equation}{section}

\begin{document}

\title[Disjointness and discrete elements in partially ordered spaces]{A note on disjointness and discrete elements in partially ordered vector spaces}

\thispagestyle{plain}

\author{Jani Jokela}
\address[J. Jokela]{Mathematics Research Centre, Tampere University, P.O. Box 692, 33101 Tampere, Finland.}
\email{jani.jokela@tuni.fi}

\begin{abstract}
The notions of disjointness and discrete elements play a prominent role in the classical theory of vector lattices. 
There are at least three different generalizations of the notion of disjointness to a larger class of partially ordered vector spaces. In recent years, one of these generalizations has been widely studied in the context of pre-Riesz spaces. The notion of $D$-disjointness is the most general of the three disjointness concepts.   
In this paper we study $D$-disjointness and the related concept of a $D$-discrete element. We establish some basic properties of $D$-discrete elements in Archimedean partially ordered spaces, and we investigate their relationship to discrete elements in the theory of pre-Riesz spaces. We then apply our results to establish the equivalence of pervasiveness and weak pervasiveness in finite-dimensional Archimedean pre-Riesz spaces.
\end{abstract}

\date{\today}

\maketitle

\section{Introduction}

We first recall some fundamental concepts. A real vector space $V$ together with a partial ordering $\leq$ is called a \emph{partially ordered vector space}, if  
\begin{equation*}\label{pospace}
u\leq v \; \implies \; u+w\leq  v+w  \, \quad \textrm{ and } \quad \,
u\leq v \; \implies \; a u\leq a v
\end{equation*}
holds for all $u,v,w\in V$ and $a \in \R_+$. A partially ordered vector space $V$ is called a \emph{vector lattice}, or a \emph{Riesz space}, if $x\vee y=\sup \{x,y\}$ and $x\wedge y=\inf\{x,y\}$ exist in $V$ for all $x,y\in V$. 
A subset $C$ of a vector space $V$ is called a \emph{cone} if \, (i)\;$\alpha C \subseteq C$ for all $\alpha \in \R_+$, \, (ii)\;\;$C+C \subseteq C$ and \, (iii) \,$C\cap (-C)=\{0\}$. For $x,y\in V$ and any cone $C$ in $V$ there is an associated partial ordering defined by $x\leq y$ if $y-x\in C$. Then $C$ is called the \emph{positive cone} for the partial ordering $\leq$. In this paper, the positive cone of a partially ordered vector space $V$ is denoted by $V_p$. 

A partially ordered space $V$ is called \emph{directed} if the positive cone $V_p$ is generating, or equivalently, for any $x,y\in V$ there exists an element $z\in V$ such that $x\leq z$ and $y\leq z$. Moreover, $V$ is called \emph{Archimedean} if $nx\leq y$ for all $n\in \N$ implies that $x\leq 0$. If $S\subset V$ then an element $x\in S$ is called a \emph{maximal element} of $S$ if $y\in S$ and $x\leq y$ implies $y=x$. A \emph{minimal element} is defined similarly.

A partially ordered vector space $V$ is called a \emph{pre-Riesz space} if for all $x,y,z\in V$ the inclusion $\{x+z,y+z\}^U\subseteq \{x,y\}^U$ implies that $z\geq 0$. 
Every pre-Riesz space is directed and every directed Archimedean partially ordered vector space is pre-Riesz \cite[Proposition 2.2.3]{pre}.

Two elements $x$ and $y$ in a vector lattice $V$ are called \emph{disjoint} if $|x|\wedge |y|=0$. Disjointness is a crucial concept in the structure theory of vector lattices. 
If we wish to extend this notion to more general partially ordered spaces, the disjointness condition must be replaced by some other, more general condition. 
There is a generalization of the notion of disjointness, where the elements $x$ and $y$ are called \emph{disjoint} if the sets $\{x+y,-(x+y)\}$ and $\{x-y,y-x\}$ have the same upper bounds. This definition of disjointness has been used in the study of pre-Riesz spaces \cite{pre}. 
In this paper we study a more general notion of \emph{$D$-disjointness}, which was first introduced in \cite{poly}. Its relationship to the other formulations of disjointness has been studied in \cite{glu, pre} with some interesting applications.  

The closely related notion of a \emph{discrete element} plays an important role in the theory of Archimedean vector lattices. This concept has also been extended to a more general setting of pre-Riesz spaces \cite{mali}. 
In Section \ref{sec:three} we introduce another type of generalization based on the notion of $D$-disjoint elements. 
We then apply these ideas in Section \ref{sec:four} to study finite-dimensional pre-Riesz spaces. It is known that a finite-dimensional Archimedean pervasive pre-Riesz space is a vector lattice \cite[Theorem 39]{mali}. We will improve on this result by showing that a finite-dimensional weakly pervasive pre-Riesz space is already a vector lattice, under the Archimedean hypothesis. This result has been proved also in \cite[Corollary 4.14]{glu} by a very different approach using topological methods, whereas the proof given in this paper is purely algebraic. 

In the existing literature, many theorems on pre-Riesz spaces are proved by embedding the pre-Riesz space into its vector lattice cover. In contrast, the present paper avoids the use of the vector lattice cover. All our results are proved using only the intrinsic properties of the positive cone of the space we are dealing with.

\section{Disjointness in partially ordered spaces}\label{sec:two}

Throughout this paper, we will use the following notation. The set of all upper bounds of a set $A$ is denoted by $A^U$, and the set of lower bounds of $A$ is denoted by $A^L$. 
The next definition contains two different notions of disjointness.  

\begin{defn}\label{def_disj}
Let $(V,\leq)$ be a partially ordered vector space.
\begin{enumerate}[(i)]
\item 
The elements $x$ and $y$ are called \emph{disjoint} 
if the sets $\{x+y,-(x+y)\}$ and $\{x-y,y-x\}$ have the same nonempty sets of upper bounds. 
This is denoted by $x \perp y$.
\item 
The elements $x,y\in V_p$ are called \emph{$D$-disjoint} if $[0,x]\cap[0,y]=\{0\}$.  
This is denoted by $x\pperp y$.
\end{enumerate} 
\end{defn}

The notion of $D$-disjointness was introduced in \cite{poly}. Disjointness was first discussed in \cite{dis}, and since then it has turned out to be a very fruitful concept in the study of pre-Riesz spaces.

It is worthwhile to make a few remarks concerning these definitions. 
The above definition of disjointness is based on the fact that two elements $x$ and $y$ in a vector lattice are disjoint if and only if $|x-y|=|x+y|$. However, this condition is not the actual definition of disjointness. By definition, two elements $x$ and $y$ in a vector lattice $V$ are disjoint if $|x|\wedge |y|=0$. If $x,y\geq 0$ this condition can be expressed as 
$$
\max\{w\in V: w\leq x \, \text{ and } \, w\leq y\}=\max\{x,y\}^L = 0.
$$
In a partially ordered space the above maximum does not necessarily exist, and the most natural step to generalize this notion is to require that, instead of the maximum, the set $\{x,y\}^L$ contains \emph{maximal elements}, and that zero is one of these maximal elements. If we denote the set of maximal elements of a set $A$ by $\mmax A$, then this condition can be stated as $0\in \mmax \{x,y\}^L$. This is equivalent to the definition of $D$-disjointness. Since we will mostly apply the definition of $D$-disjointness in this form, it is convenient to state this simple observation as a lemma. 

\begin{lemma}\label{d_disj_def}
The elements $x,y\in V_p$ are $D$-disjoint if and only if $0\in \mmax \{x,y\}^L$. 
\end{lemma}

\begin{bproof}
We note that $[0,x]\cap[0,y]=\{w\geq 0: w\leq x , \, w\leq y\}=V_p\cap \{x,y\}^L$. Hence $x$ and $y$ are $D$-disjoint if and only if $V_p\cap \{x,y\}^L=\{0\}$. Now if $x$ and $y$ are $D$-disjoint, and $w\in \{x,y\}^L$ then $w\geq 0$ implies $w=0$, so that $0\in\mmax \{x,y\}^L$. Conversely, if $0\in\mmax \{x,y\}^L$ then $0$ is the only positive element in $\{x,y\}^L$. Thus $V_p\cap \{x,y\}^L=\{0\}$ and $x$ and $y$ are $D$-disjoint.
\end{bproof}

In general, $D$-disjointness is different than the notion of disjointness, although under certain conditions 
these two notions coincide. 
In fact, $D$-disjointness always implies disjointness, and hence $D$-disjointness is a more general condition than disjointness. This result has been proved in the pre-Riesz space setting in \cite[Proposition 4.1.9]{pre}, where it was also shown that these two notions are equivalent in a pre-Riesz space with the Riesz decomposition property. 
This holds, in particular, in a vector lattice.

It is also worth pointing out that another disjointness-type concept was discussed in \cite{glu}, where the following condition was considered: 
\begin{equation}\label{diseq}
[-x,x]\cap[-y,y]=\{0\}.
\end{equation}
It was proved in \cite[Proposition 4.16]{glu} that $x\perp y$ implies condition \eqref{diseq}, which in turn implies that $x\pperp y$. Furthermore, it was shown by way of examples that none of these implications can be reversed. Hence the condition \eqref{diseq} is intermediate between disjointness and $D$-disjointness.

The existence of maximal elements of the set $\{x,y\}^L$ depends on the properties of the positive cone. It is therefore of some interest to consider those positive cones for which the existence of maximal elements is guaranteed.

\begin{defn}
Let $(V,\leq)$  be a partially ordered vector space. The positive cone $V_p$ is called an \emph{$M$-cone} if for all $x,y\in V$ the set $\{x,y\}^L$ contains maximal elements. 
The set of maximal elements of $\{x,y\}^L$ is denoted by $\mmax \{x,y\}^L$. 
\end{defn}

We observe that $V_p$ is an $M$-cone if and only if for all $x,y\in V$ the set defined by $\{x,y\}^U=\{w\in V: w\geq x, \; w\geq y\}$ 
contains minimal elements. 
This follows immediately from the fact that if $m$ is a maximal element of a set $A$ then $-m$ is a minimal element of the set $-A$. The set of minimal elements of the set $\{x,y\}^U$ is denoted by $\mmin \{x,y\}^U$.

The sets of maximal and minimal elements have the following translation properties, which are easily verified.
\begin{equation}\label{transl1}
\mmax \{x+z,y+z\}^L = z+ \mmax \{x,y\}^L
\end{equation}
\begin{equation}\label{transl2}
\mmin \{x+z,y+z\}^U = z+ \mmin \{x,y\}^U
\end{equation}
Moreover, if $x,y\in V_p$ and $0\leq \alpha \in \R$ then 
\begin{equation}\label{scalar_multipl}
w\in \mmax \{\alpha x,\alpha y\}^L \quad \text{if and only if} \quad \frac{1}{\alpha}w \in \mmax \{x,y\}^L. 
\end{equation}
We have the following simple characterization for $M$-cones.

\begin{prop}\label{joku}
Let $V$ be a directed partially ordered space. Then the following are equivalent.
\begin{enumerate}[(a)]
\item
The positive cone $V_p$ is an $M$-cone.
\item
The set $\{x,y\}^L$ contains maximal elements for all $x,y\in V_p$. 
\end{enumerate}
\end{prop}

\begin{bproof}
Clearly (a) implies (b). Assume that (b) holds and let $x,y\in V$. Since $V$ is directed, there exists an element $z\in V$ such that $-x\leq z$ and $-y\leq z$. Then $z+x\geq 0$ and $z+y\geq 0$, and by our hypothesis, the set $\{z+x,z+y\}^L$ has maximal elements. If $w\in \mmax \{z+x,z+y\}^L$ then $w-z \in \mmax \{x,y\}^L$, thus proving that $V_p$ is an $M$-cone. 
\end{bproof}

There are partially ordered spaces in which nonzero $D$-disjoint elements do not exist (c.f. Example \ref{not_mcone}), but the situation is different if $V_p$ is an $M$-cone. 

\begin{prop}
Let $V$ be a partially ordered space such that $V_p$ is an  
$M$-cone. If $\dim V \geq 2$ then $V$ contains nontrivial $D$-disjoint elements. 
\end{prop}

\begin{bproof}
Let $x,y\in V$ be nonzero elements such that $x\neq cy$ ($c\in\R$). 
Let $z\in \mmax \{x,y\}^L$. If $u=x-z$ and $v=y-z$  then $0\in \mmax \{u,v\}^L$, by \eqref{transl1}. Hence $(x-z)\pperp (y-z)$, by Lemma \ref{d_disj_def}. 
\end{bproof}

The converse of the above does not hold. That is, $V_p$ may contain $D$-disjoint elements even if $V_p$ is not an $M$-cone.

\begin{example}
Let $V=\R^3$ with the positive cone $V_p=\{(x,y,z): x>0, y\geq 0, z\geq 0\}\cup \{(0,y,z): y\geq 0, z\geq 0\}$. Then $V_p$ is not an $M$-cone. For example, if $u=(1,1,1)$ and $v=(0,2,0)$ then the set $\{u,v\}^L$ is nonempty but it has no maximal elements. 
On the other hand, if $a=(0,1,0)$ and $b=(0,0,1)$ then $a\pperp b$. 
Note also that $V$ is not pre-Riesz. 
To see this, let $x=(0,1,0)$, $y=(1,1,0)$ and $z=(1,0,0)$. If $w=(w_1,w_2,w_3)\in \{x+z,y+z\}^U$ then $w\geq (1,1,0)$ and $w\geq (2,1,0)$ so that $w_1> 2$ and $w_2\geq 1$. 
Then $w-x=(w_1,w_2-1,w_3)\in V_p$ and $w-y=(w_1-1,w_2-1,w_3)\in V_p$, so $w\in \{x,y\}^U$ but $z\notin V_p$.  
\end{example}

Although the definition of $D$-disjointness is limited to positive elements, it is more general than the notion of disjointness, as far as positive elements are concerned. 
The fact that $D$-disjointness implies disjointness was proved for pre-Riesz spaces in \cite[Proposition 4.1.9]{pre}, where the proof used methods involving embedding the pre-Riesz space into its vector lattice cover. We will show below in Theorem \ref{alternative} that this result holds in any partially ordered space using elementary methods. For the proof we need the following simple lemma.

\begin{lemma}\label{disjointnesslemma1}
Let $(V,\leq)$  be a partially ordered vector space. 
Then $x\pperp y$ if and only if $\alpha x\pperp \alpha y$ for all $0\leq \alpha \in \R$.
\end{lemma}

\begin{bproof} 
Assume that $x\pperp y$ and $0\leq w\in \{\alpha x , \alpha y\}^L$. The case $\alpha =0$ is trivial, so let $\alpha >0$. Then $0\leq \frac{1}{\alpha}w\in \{x , y\}^L$, by \eqref{scalar_multipl}. By assumption, $0\in \mmax \{x , y\}^L$ which implies that $\frac{1}{\alpha}w\leq 0$. Consequently, $w=0$ and hence $0\in \mmax \{\alpha x , \alpha y\}^L$. Thus $\alpha x\pperp \alpha y$. The converse implication is evident by choosing $\alpha =1$. 
\end{bproof}

\begin{thm}\label{alternative}
Let $V$ be a partially ordered vector space and $x,y\in V_p$. 
Then $x\perp y$ implies $x\pperp y$. 
\end{thm} 

\begin{bproof}
Let $A=\{x+y,-x-y\}$ and $B=\{x-y,y-x\}$. First we note that $B=-B$ and thus $w\in \mmin B^U$ if and only if $-w\in\mmax B^L$. This in turn is equivalent to $0\in \mmax \{x-y+w,y-x+w\}^L$, or $(x-y+w)\pperp (y-x+w)$. 
Hence $w\in \mmin B^U$ if and only if $(x-y+w)\pperp (y-x+w)$. 

Assume then that $x\perp y$. Since $x,y\geq 0$, we have $x+y=\min A^U$. Since $A^U=B^U$, it follows that 
$w=x+y=\min B^U$. We showed above that this is equivalent to $(x-y+w)\pperp (y-x+w)$. But $x-y+w=2x$ and $y-x+w=2y$, and so $2x\pperp 2y$. By Lemma \ref{disjointnesslemma1} this is equivalent to $x\pperp y$, and the proof is complete. 
\end{bproof}

The converse of the preceding theorem does not hold in general, as shown in \cite[Example 4.1.8]{pre}. 
However, the converse implication in Theorem \ref{alternative} holds if $V$ has the Riesz decomposition property \cite[Proposition 4.1.9 (ii)]{pre}. 
(Note that in \cite{pre} this result was stated for pre-Riesz spaces but their proof makes no use of the assumption that $V$ is pre-Riesz.)

\section{Discrete elements}\label{sec:three}

We recall that an element $x$ in a vector lattice is called \emph{discrete} if $0\leq a\leq |x|$ and $0\leq b\leq |x|$ with $|a|\wedge |b|=0$ implies that $a=0$ or $b=0$. In other words, $x$ is discrete if there are no nonzero disjoint elements $a$ and $b$ satisfying $0\leq a\leq |x|$ and $0\leq b\leq |x|$. A generalization of a discrete element to pre-Riesz spaces was given in \cite{mali}, where the lattice-disjointness was replaced by the generalized disjointness condition of Definition \ref{def_disj}. The notion of $D$-disjointness now enables us to give a corresponding definition of a discrete element.

\begin{defn}\label{def_atom}
Let $(V,\leq)$  be a partially ordered space and $x,y\in V_p$. 
\begin{enumerate}[(i)]
\item
An element $x\neq 0$ is called an \emph{atom}  
if $0\leq a\leq x$ implies that $a=\alpha x$ for some $\alpha \geq 0$
\item
An element $x\neq 0$ is called \emph{discrete}  
if there are no nonzero disjoint elements $a$ and $b$ such that $0\leq a\leq x$ and $0\leq b\leq x$. 
\item
An element $x\neq 0$ is called \emph{$D$-discrete}  
if there are no nonzero $D$-disjoint elements $a$ and $b$ such that $0\leq a\leq x$ and $0\leq b\leq x$.
\end{enumerate}
\end{defn}

We have excluded the zero element from the above definition for convenience, although zero trivially satisfies the conditions of Definition \ref{def_atom}. 

The above definition of an atom is a well-established concept in the theory of partially ordered spaces. 
The generalized notion of a discrete element was introduced in \cite{mali}. To the best of the author's knowledge, the concept of a $D$-discrete element has not been studied in the existing literature.

It follows immediately from Theorem \ref{alternative} that every $D$-discrete element is discrete, but the converse is not true (c.f. Example \ref{ex_discrete}). 
On the other hand, we have the following relationship between atoms and $D$-discrete elements.

\begin{prop}\label{atom_discrete}
Let $V$ be a partially ordered space. Then every atom is $D$-discrete.
\end{prop}

\begin{bproof}
Let $x$ be an atom and let $a,b$ be $D$-disjoint elements such that   
$0\leq a \leq x$ and $0\leq b \leq x$. Then $a=\alpha x$ and $b=\beta x$. Thus if $c=\min\{\alpha,\beta\}$ then $cx\geq 0$, and so $cx\in V_p\cap \{a,b\}^L =\{0\}$.   
Hence $c=0$, and consequently, $a=0$ or $b=0$. Thus $x$ is $D$-discrete.
\end{bproof}

The converse of the preceding result holds if $V$ is Archimedean and $V_p$ is an $M$-cone.

\begin{thm}\label{discrete_atom}
Let $V$ be an Archimedean 
partially ordered space such that $V_p$ is an $M$-cone. Then an element $x\in V$ is an atom if and only if $x$ is a $D$-discrete element. 
\end{thm}

\begin{bproof}
By Proposition \ref{atom_discrete} it remains to show that if $x$ is $D$-discrete then $0\leq a \leq x$ implies that $a=\alpha x$ for some $\alpha \geq 0$. The case $a=0$ is trivial so let us assume that $a\neq 0$. 

Since $V$ is Archimedean, the set $S=\{c\in \R: ca\leq x\}$ is bounded, so that $\alpha=\sup S>0$ exists and $\alpha a\leq x$. Let $v=x-\alpha a \geq 0$ and let $\lambda\in\R$ be an arbitrary number such that $0<\lambda\leq 1$. Then since $a\neq 0$, we have $\lambda a \neq 0$ and the 
inequalities $0\leq \lambda a \leq a \leq x$ and $0\leq v \leq x$ hold. 
  
Assume then that $z\in \mmax \{v,\lambda a\}^L$. Then $z\geq 0$ and $0\in \mmax \{v-z,\lambda a-z\}^L$, and so $(v-z)\pperp (\lambda a -z)$. 
We also have $v-z\leq x$ and $\lambda a-z\leq x$, and since $x$ is $D$-discrete, it follows that $\lambda a -z=0$ or $v-z=0$.

If $z=\lambda a$ then $\lambda a \leq v=x-\alpha a$, and so $(\alpha +\lambda)a\leq x$. But this contradicts the definition of $\alpha$, so we must have $z\neq \lambda a$, and hence $v-z=0$. But then $v=z< \lambda a$ for all $\lambda\in (0,1]$. Since $V$ is Archimedean, this implies that $v=0$, or $x=\alpha a$. 
\end{bproof}

\begin{cor}\label{discrete_cor}
Let $V$ be an Archimedean partially ordered space such that $V_p$ is an $M$-cone. If $x$ and $y$ are $D$-discrete elements then either $x$ and $y$ are $D$-disjoint or $x=\alpha y$ for some $\alpha\geq 0$.
\end{cor}

\begin{bproof}
If $x$ and $y$ are not $D$-disjoint then there exists an element 
$z$ such that $0<z\leq x$ and $0<z\leq y$. Since $x$ and $y$ are atoms, there exists $c,d >0$ such that $x=c z$ and $y=d z$. Then $x= \frac{c}{d} y$.
\end{bproof}

The preceding corollary can be reformulated as follows. 

\begin{cor}\label{discrete_cor2}
Let $V$ be an Archimedean partially ordered space such that $V_p$ is an $M$-cone. Then $x$ and $y$ are $D$-disjoint atoms if and only if $x$ and $y$ are linearly independent.
\end{cor}

%
%

\begin{example}\label{not_mcone}   
Let $V=\R^2$ with the positive cone $V_p=\{(x,y): x>0, \; y>0\}\cup \{(0,0)\}$. Then every positive element of $V$ is $D$-discrete (because $V$ does not contain any $D$-disjoint elements that are both nonzero) but $V$ contains no atoms. Notice that $V$ is not Archimedean and $V_p$ is not an $M$-cone. 
\end{example}

Theorem \ref{discrete_atom} and Corollary \ref{discrete_cor} are generalizations of corresponding results in Archimedean vector lattices (see \cite[\S 13]{vul}). However, it is worth pointing out that these results do not hold for discrete elements in more general Archimedean spaces. 
The following example clarifies the difference between the notions of discrete and $D$-discrete elements. 
In particular, it shows that discrete elements are not necessarily $D$-discrete.

\begin{example}\label{ex_discrete}
Consider an Archimedean partially ordered space $V=\R^3$ with the positive cone $K=\{c_1v_1+c_2v_2+c_3v_3+c_4v_4: c_k\geq 0\}$ where $v_1=(1,0,1)$, $v_2=(0,1,1)$, $v_3=(-1,0,1)$ and $v_4=(0,-1,1)$. It was shown in \cite[Example 31]{mali} that the atoms in $V$ are precisely the vectors $cv_k$ for $c>0$. It was also shown that for the vector $x=v_1+v_2$ there are no nonzero elements $u,v$ such that $0\leq u\leq x$, $0\leq v\leq x$ and $u\perp v$. Hence $x$ is discrete but not an atom. However, $v_1$ and $v_2$ are nonzero $D$-disjoint elements (see \cite[Example 4.8]{glu}) such that $0\leq v_1\leq x$ and $0\leq v_2\leq x$.  
Thus $x$ is not $D$-discrete. 
\end{example}

We note that if $V_p$ is an $M$ cone then $V_p$ is generating, that is, $V$ is a directed space. Every lattice cone is obviously an $M$-cone. 
It is known that every directed Archimedean partially ordered space is pre-Riesz \cite[Proposition 2.2.3]{pre}. 
Hence every Archimedean partially ordered space whose positive cone is an $M$-cone is pre-Riesz. In the finite-dimensional case we have the following characterization.

\begin{thm}\label{preriesz_char1}
A finite-dimensional Archimedean partially ordered vector space $V$ is a pre-Riesz space if and only if $V_p$ is an $M$-cone. 
\end{thm}

\begin{bproof}
By the preceding remark it is sufficient to prove the other implication. 
Let $V$ be pre-Riesz and $x,y\in V$. Then $V$ is directed so that $\{x,y\}^U$ is non-empty. Let $z\in\{x,y\}^U$ and define $S=[x,z]\cap [y,z]$. Since $V$ is finite-dimensional and Archimedean it follows that if $V$ is equipped with a norm then the positive cone is closed and it contains interior points. Consequently, $V$ admits a strictly positive linear functional $\phi$. Moreover, since $V_p$ is closed, every order interval in $V$ is compact so that the set $S$ is also compact. Thus $\phi$ attains its minimum at some $u\in S$. Then $u$ is in fact a minimal element of $\{x,y\}^U$. Indeed, if $v\in \{x,y\}^U$ is such that $v\leq u$, then $v\in S$. Now if $v\neq u$ then $u-v>0$, and consequently, $\phi(u-v)>0$. But then $\phi(v)<\phi(u)$, a contradiction. Hence $u\in \mmin \{x,y\}^U$ and so $V_p$ is an $M$-cone. 
\end{bproof}


\begin{remark}
There are non-Archimedean finite-dimensional pre-Riesz \linebreak spaces whose positive cone is an $M$-cone. For example, if $V=\R^2$ with the lexicographic ordering then $V$ is a totally ordered vector lattice, so that the positive cone is certainly an $M$-cone, but $V$ is not Archimedean. On the other hand, finite-dimensionality is essential in Theorem \ref{preriesz_char1}. For example, let $V$ be the space of all differentiable functions on $[0,1]$. If $f(x)=x$ and $g(x)=1-x$ then the set $\{f,g\}^U$ does not contain minimal elements, and so $V_p$ is not an $M$-cone. 
\end{remark}

\section{Atoms in a weakly pervasive pre-Riesz space}\label{sec:four}

Due to Theorem \ref{preriesz_char1}, most results of this paper are applicable to all finite-dimensional Archimedean pre-Riesz spaces. 
For the applications that follow, we recall that 
a pre-Riesz space $V$ is called \emph{weakly pervasive} if any two $D$-disjoint elements in $V_p$ are disjoint. Consequently, an element $x$ in a weakly pervasive pre-Riesz space is $D$-discrete if and only if $x$ is discrete. Therefore, in the context of pre-Riesz space theory, Theorem \ref{discrete_atom} and Corollary \ref{discrete_cor2} have the following immediate consequences. \\

\begin{cor}\label{last_cor}
Let $V$ be a finite-dimensional Archimedean weakly pervasive pre-Riesz space. Then the following hold. 
\begin{enumerate}[(a)]
\item
An element $x\in V$ is an atom (i.e. $D$-discrete) if and only if $x$ is discrete. 
\item
Elements $x$ and $y$ are disjoint atoms if and only if $x$ and $y$ are linearly independent.
\item
Any set of pairwise linearly independent atoms is linearly independent.
\end{enumerate}
\end{cor}

\begin{bproof}
By the preceding remarks, the only item that requires proof is (c). Let $S=\{x_1, x_2,\ldots ,x_k\}$ be a pairwise linearly independent set of atoms. Suppose that $S$ is linearly dependent. Then we may assume that $x_1=\sum_{i=2}^k c_i x_i$, where at least one of the coefficients, say $c_j$, is nonzero. Then 
$$
0\leq c_jx_j \leq \sum_{i=2}^k c_i x_i =x_1.
$$
But $c_j x_j$ and $x_1$ are both atoms, so by Corollaries \ref{discrete_cor} and \ref{discrete_cor2} we have $c_jx_j =cx_1$ for some $c>0$. This contradicts the assumption that $S$ is pairwise linearly independent. Hence $S$ is linearly independent.
\end{bproof}

The statements of Corollary \ref{last_cor} were proved in \cite[Theorem 32 and Propositions 35 and 36]{mali} for the class of Archimedean pervasive pre-Riesz spaces. Note that Corollary \ref{last_cor}(c) does not hold if $V$ is not weakly pervasive. This is shown by Example \ref{ex_discrete}, where the set of atoms $\{v_1, v_2, v_3, v_4\}$ is pairwise linearly independent but not linearly independent. 

For the proof of our next result we recall that a sub-cone $C$ of $V_{p}$ is called \emph{order convex} if $0\leq w \leq z$ with $z\in C$ implies $w\in C$. An order-convex sub-cone is also called a \emph{face}. 
The face $C_x$ of $V_{p}$ generated by an element $x\geq 0$ is given by 
$C_x=\{w: 0\leq w\leq cx \; \text{ for some } c\geq 0\}$.

It is not immediately clear whether or not $D$-discrete elements always exist but 
if $V$ is Archimedean and finite dimensional then we have the following result.

\begin{thm}\label{discrete_exist1}
Let $(V,\leq)$ be a finite dimensional Archimedean pre-Riesz space. 
If $F_1$ is a face of $V_p$ then for any nonzero element $x\in F_1$ there exists an atom 
$z$ such that $0\leq z\leq x$. In particular, this holds if $F_1=V_p$. 
\end{thm}

\begin{bproof}
Let $F_1$ be a face of $V_{p}$ (possibly $F_1=V_{p}$) 
and $x_1\in F_1$ a nonzero element. Suppose that $a_1$ is a nonzero element such that $0\leq a_1 \leq x_1$. If every such element $a_1$ is a scalar multiple of $x_1$ then by Theorem \ref{discrete_atom} $x_1$ is an atom, and we are done. 
If $a_1$ is not a scalar multiple of $x_1$ then we define $x_2=x_1-\alpha_1 a_1$ where $\alpha_1=\sup\{c\in \R: ca_1\leq x_1\}$ (Note that such number $\alpha_1$ exists since $V$ is Archimedean). 
Then $x_2$ is nonzero with $0\leq x_2\leq x_1$, and if $F_2$ is the face generated by $x_2$ then $x_1\notin F_2$. For otherwise we would find some $c>1$ such that 
$x_1\leq c(x_1-\alpha_1 a_1)$, or equivalently $\alpha_1 a_1 \leq \frac{c-1}{c}x_1 \leq x_1$, which contradicts the definition of $\alpha_1$ since $\frac{c-1}{c}<1$. Hence $F_2\subset F_1$ and $F_2\neq F_1$, so that $F_2$ is a proper face of $F_1$.  

Next, assume that $a_2$ is a nonzero element such that $0\leq a_2 \leq x_2$. If $a_2$ is not a scalar multiple of $x_2$ then we define $x_3=x_2-\alpha_2 a_2$ where $\alpha_2=\sup\{c\in \R: ca_2\leq x_2\}$. Similar reasoning as above shows that $0\neq x_3\leq x_2$ and the face $F_3$ generated by $x_3$ is a proper face of $F_2$, that is, $F_3\subset F_2\subset F_1$. In particular, $x_2\notin F_3$. 

Continuing in this manner, we produce elements $x_1\geq x_2 \geq \ldots \geq x_k\neq 0$ and a corresponding sequence of faces $F_k\subset F_{k-1} \subset \ldots \subset F_1$ such that each $F_i$ is a proper face of $F_{i-1}$, and $x_i\notin F_{i+1}$ for all $i=1,\ldots, k-1$. 

Now the set $\{x_1,x_2,\ldots x_k\}$ is linearly independent. To see this, assume that $\sum_{i=1}^k c_i x_i =0$ such that some of the coefficients are nonzero. Let $m$ be the smallest index such that $c_m\neq 0$. If $m=k$ then $x_k=0$, a contradiction. If $m<k$ we have $x_m=\sum_{i=m+1}^k c_i x_i$. But then $x_m\in F_{m+1}$ which is again a contradiction. Hence $\{x_1,x_2,\ldots x_k\}$ is linearly independent.

It follows that if the dimension of $V$ is finite, say $n$, then the above process must eventually terminate, and we find that for some $k\leq n$, if $a_k$ is an element such that $0\leq a_k \leq x_k$ then we must have $a_k=\alpha x_k$ for some $\alpha \geq 0$. Hence we have found an atom $x_k$ such that $0\leq x_k\leq x_1$. 
\end{bproof}

We shall call two atoms $x$ and $y$ \emph{essentially distinct} if $x\neq cy$ for all $c>0$. 
With the preceding theorem, we can now prove the following.

\begin{thm}\label{main1}
Let $V$ be an $n$-dimensional Archimedean weakly pervasive pre-Riesz space. 
Then $V$ contains precisely $n$ essentially distinct atoms, and these atoms form a basis of $V$. Moreover, if $\{x_1, \ldots , x_k\}$ is such basis then $V_p=\textrm{pos}\{x_1, \ldots , x_k\}$. In fact, every $y\in V_p$ can be written as $y=\sum_{i=1}^{n} \alpha_i x_i$, where
\begin{equation}\label{esitys}
\alpha_j=\sup\{c\in \R: cx_j \leq y- \sum_{i=1}^{j-1} \alpha_i x_i\}.
\end{equation} 
\end{thm}

\begin{bproof}
Let $\dim (V)=n$ and let $S=\{x_1, \ldots , x_k\}$ be a maximal linearly independent set of atoms. By Theorem \ref{discrete_exist1} the set $S$ is nonempty. Clearly we have $k\leq n$. Suppose that $k<n$. Then there exists $y\in V_p$ such that $y\notin \Span S$. We now use a similar method as in the proof of Theorem \ref{discrete_exist1}. Define the numbers $\alpha_j$ according to \eqref{esitys}. 
Then $z= y-\sum_{i=1}^{k} \alpha_i x_i$ is a positive nonzero element, and by Theorem \ref{discrete_exist1} there exists an atom $x$ such that $x\leq z$. 
Since each $x_j$ is an atom, it follows by Corollary \ref{discrete_cor} that for each $x_j$ we have either $x\pperp x_j$ or $x=c_jx_j$ for some $c_j>0$. Now if $x\pperp x_j$ holds for all $j=1,\ldots, k$ then by Corollary \ref{discrete_cor2}, the set $S_1=\{x_1, \ldots , x_k, x\}$ is pairwise linearly independent, and so by Corollary \ref{last_cor}(c) the set $S_1$ is linearly independent. But this contradicts the maximality of $S$. Hence we must have $x=c_jx_j$ for some $j\in\{1,\ldots, k\}$ and $c_j>0$. But then the inequality $x\leq z$ becomes
$$
0\leq c_j x_j \leq y- \sum_{i=1}^{k} \alpha_i x_i \leq y -\sum_{i=1}^{j} \alpha_i x_i,
$$
and this implies that
$$
(\alpha_j + c_j)x_j \leq y- \sum_{i=1}^{j-1} \alpha_i x_i.
$$
This contradicts the definition of $\alpha_j$. Hence we must have $k=n$, and so $S$ is a basis of $V$ consisting of essentially distinct atoms. 
This also shows that for any $y\in V_p$ the element $z= y-\sum_{i=1}^{n} \alpha_i x_i$ must be zero, for if $z$ were nonzero then we could find an atom $x$ such that $x\leq z$. Then repeating the preceding argument yields a contradiction. Hence $y=\sum_{i=1}^{n} \alpha_i x_i$. Clearly, $\alpha_j \geq 0$ for all $j=1, \ldots , n$ and so $V_p=\textrm{pos}\{x_1, \ldots , x_k\}$. Finally, since any set of essentially distinct atoms is pairwise $D$-disjoint, and hence pairwise linearly independent, it follows by Corollary \ref{last_cor}(c) that the set is linearly independent and thus it cannot contain more than $n$ atoms. 
\end{bproof}

A pre-Riesz space $V$ is called \emph{pervasive} if for every $x\in V$ such that $x\nleq 0$ there exists $y\in V_p \setminus\{0\}$ such that every positive upper bound of $x$ is also an upper bound of $y$. 
The preceding theorem implies that in finite-dimensional Archimedean pre-Riesz spaces the notions of pervasiveness and weak pervasiveness coincide. Indeed, we have the following characterization.

\begin{cor}\label{main2}
Let $V$ be an $n$-dimensional Archimedean directed partially ordered vector space. Then the following statements are equivalent.
\begin{enumerate}[(a)]
\item
There exists a basis $\{x_1, \ldots , x_n\}$ of $V$ such that $V_p=\textrm{pos}\{x_1, \ldots , x_n\}$.
\item
$V$ is a vector lattice.
\item
$V$ has the Riesz decomposition property.
\item
$V$ is a pervasive pre-Riesz space.
\item
$V$ is a weakly pervasive pre-Riesz space.
\end{enumerate}
\end{cor}

\begin{bproof}
The equivalence of statements (a), (b) and (c) holds by \cite[Theorem 1.7.8]{pre}. Statements (b) and (d) are equivalent by \cite[Theorem 39]{mali}. The implication (b)$\implies$(e) follows from \cite[Theorem 4.1.9]{pre}, and (e)$\implies$(a) was proved above in Theorem \ref{main1}.
%
\end{bproof}

\begin{remark}
The implication (e)$\implies$(b) of the preceding corollary was proved in \cite[Corollary 4.14]{glu} using very different techniques from the theory of ordered Banach spaces.
\end{remark}

\bibliographystyle{plain}

\end{document}